
\input amssym.def
\input amssym
\magnification=1100
\baselineskip = 0.23truein
\lineskiplimit = 0.01truein
\lineskip = 0.01truein
\vsize = 8.5truein
\voffset = 0.2truein
\parskip = 0.10truein
\parindent = 0.3truein
\settabs 12 \columns
\hsize = 5.8truein
\hoffset = 0.2truein

\setbox\strutbox=\hbox{%
\vrule height .708\baselineskip
depth .292\baselineskip
width 0pt}
\font\caps=cmcsc10

\font\bigtenrm=cmr10 at 14pt

\def\sqr#1#2{{\vcenter{\vbox{\hrule height.#2pt
\hbox{\vrule width.#2pt height#1pt \kern#1pt
\vrule width.#2pt}
\hrule height.#2pt}}}}
\def\square{\mathchoice\sqr46\sqr46\sqr{3.1}6\sqr{2.3}4}

\centerline{\bigtenrm EXPANDERS, RANK AND GRAPHS OF GROUPS}

\tenrm
\vskip 14pt
\centerline{MARC LACKENBY}
\vskip 18pt

\centerline{\caps 1. Introduction}
\vskip 6pt

A central principle of this paper is that, for a
finitely presented group $G$, the algebraic properties
of its finite index subgroups should be reflected by
the geometry of its finite quotients. These quotients
can indeed be viewed as geometric objects, in
the following way. If we
pick a finite set of generators for $G$, these map to
a generating set for any finite quotient and hence
endow this quotient with a word metric. This metric
of course depends on the choice of generators, but
if we were to pick another set of generators for $G$,
the metrics on the quotients would change by a
bounded factor. Thus, although the metric on any
given finite quotient is unlikely to be useful, the metrics
on the whole collection of finite quotients have a good
deal of significance.

The above principle is inspired by a common theme in
manifold theory: that the geometry and topology of a Riemannian
manifold should relate to its fundamental group. Thus,
if we pick a compact Riemannian manifold with fundamental
group $G$ (which is possible since $G$ is finitely presented), 
and let $M_i$ be the covering space corresponding
to a finite index normal subgroup $G_i$, then the geometry of $M_i$ should
have consequences for the algebraic structure of $G_i$.
But $M_i$ is coarsely approximated by the word metric on the
quotient $G/G_i$.

In this paper, we will focus on the geometric properties of the quotient
groups that relate to their Cheeger constant. Recall
that this is defined as follows. Fix a finite set $S$ of generators for
$G$, and let $X_i$ be the Cayley graph of $G/G_i$ with
respect to $S$. The {\sl Cheeger constant} of $X_i$, denoted
$h(X_i)$, is defined to be
$$\min \left\{ {|\partial A| \over |A|} : A \subset V(X_i)
\hbox{ and } 0 < |A| \leq |V(X_i)|/2 \right\}.$$
Here, $V(X_i)$ is the vertex set of $X_i$, and
$\partial A$ denotes the set of edges joining
a vertex in $A$ to one not in $A$. The group $G$
is said to have {\sl Property $(\tau)$} with respect
to a collection $\{ G_i \}$ of finite index normal
subgroups if the Cheeger constants $h(X_i)$ are bounded away from zero.
The graphs $X_i$ then form what is known as an
{\sl expanding family} or an {\sl expander}.
There is a remarkably rich theory relating to
Property $(\tau)$ [7]. Possibly its most striking aspect
is that it has so many equivalent definitions, drawing
on many different areas of mathematics, including
graph theory, differential geometry and representation
theory. It is, in general, very difficult
to construct explicit expanding families of graphs
with bounded valence ([7],[10]). A consequence
of this paper is that they are probably very common.

We will focus mainly on two algebraic properties of the
subgroups $G_i$. The first is whether or not $G_i$
decomposes as an amalgamated free product or HNN extension;
in other words, whether or not $G_i$ admits a non-trivial
decomposition as a graph of groups.
Such groups play a central r\^ole in combinatorial group
theory. The second is a new concept related to their
rank. The {\sl rank} of a group $G$ is the minimal
size of a generating set, and is denoted $d(G)$. When $G_i$ is a finite
index subgroup of $G$, the Reidermeister-Schreier
process [9] gives of a collection of $(d(G) - 1)[G:G_i] + 1$
generators for $G_i$. Thus,
$$d(G_i) - 1 \leq (d(G) - 1) [G:G_i].$$
When $\{ G_i \}$ is a collection of finite
index subgroups of $G$,
the {\sl rank gradient} of the pair $(G, \{ G_i \})$ measures
the strictness of this inequality.
It is defined to be
$$\inf_i \left \{ { d(G_i) - 1 \over [G:G_i]} \right \}.$$
The {\sl rank gradient} of $G$ is defined by taking $\{ G_i \}$
to be the set of all its finite index normal subgroups.

A collection $\{ G_i \}$ of subgroups of a group $G$ is termed a
{\sl lattice} if, whenever $G_i$ and $G_j$ are in the
collection, so is $G_i \cap G_j$.
We can now state the main theorem of this paper.

\noindent {\bf Theorem 1.1.} {\sl Let $G$ be a finitely
presented group, and let $\{ G_i \}$ be a lattice
of finite index normal subgroups. Then at least one of the
following holds:
\item{1.} $G_i$ is an amalgamated free product
or HNN extension for infinitely many $i$;
\item{2.} $G$ has Property $(\tau)$ with respect
to $\{ G_i \}$;
\item{3.} the rank gradient of $(G, \{ G_i \})$ is zero.}

We will prove this theorem in \S2.

It is conclusion (3) in Theorem 1.1 that is least familiar.
It raises the question of which groups have non-zero
rank gradient and which do not. In \S3, we will investigate
rank gradient quite thoroughly.
It seems likely that, among abstract groups
in general, those
with zero rank gradient are rather special, because
they have a sequence of finite index subgroups with 
relatively small
generating sets not arising from the Reidermeister-Schreier
process. However, there are large classes of groups
with zero rank gradient. For example, mapping
tori always have zero rank gradient. Another
class consists of the $S$-arithmetic groups
for which the congruence kernel is trivial.
A familiar example is ${\rm SL}(n, {\Bbb Z})$, for $n > 2$. 
A third class arises from
the following result, which we shall prove in \S3. This is a consequence
of a slightly strengthened version of Theorem 1.1.

\noindent {\bf Theorem 1.2.} {\sl Any finitely presented, residually finite,
amenable group has non-positive rank gradient.}

On the other hand, to prove that a given group has non-zero
rank gradient is usually rather hard, since it is often difficult
to find good lower bounds on the rank of a group.
An obvious class of examples are groups with
deficiency more than one. We shall also show that the free product of
two non-trivial groups (not both ${\Bbb Z}/ 2{\Bbb Z}$) has
non-zero rank gradient.

It should not be assumed that the three possible conclusions
in Theorem 1.1 are mutually exclusive. In \S4, we investigate
the possible combinations of these properties that
can arise.

Theorem 1.1 is a group-theoretic version of a 
result in 3-manifold theory: Theorem 1.7 of [4].
Instead of dealing with rank gradient, this used
a related notion, the Heegaard gradient of a
3-manifold, which measures the rate at which
the Heegaard genus of the manifold's finite-sheeted
covering spaces grows as a function of the
covering degree. The purpose of Theorem 1.7 in [4]
is that it represents part of a programme for
proving the virtually Haken conjecture, which 
is a key unsolved problem about 3-manifolds.
This asserts that, when $G$ is the fundamental
group of a closed hyperbolic 3-manifold, some
finite index subgroup $G_i$ should admit a
non-trivial decomposition as a graph of groups.

This is the first in a pair of papers, which use
the geometry and topology of finite Cayley graphs
as a tool in group theory. In the second [5], we prove the
following purely algebraic result.

\noindent {\bf Theorem 1.3.} {\sl Let $G$ be a finitely
presented group. Then the following are equivalent:
\item{1.} some finite index subgroup of $G$ admits
a surjective homomorphism onto a non-abelian free group;
\item{2.} there exists a sequence $G_1 \geq G_2 \geq \dots$
of finite index subgroups of $G$, each normal in $G_1$, such that
\itemitem{(i)} $G_i/G_{i+1}$ is abelian for all $i \geq 1$;
\itemitem{(ii)} $\lim_{i \rightarrow \infty} 
((\log [G_i : G_{i+1}]) / [G:G_i]) = \infty$;
\itemitem{(iii)} $\limsup_i (d(G_i/G_{i+1}) / [G:G_i])  > 0$.
}

The difficult part of this theorem is the implication
$(2) \Rightarrow (1)$. In \S5, we establish a partial
result in this direction. Using Theorem 1.1, we show
that $(2)$ implies that $G_i$ is a non-trivial
graph of groups for all sufficiently large $i$.
The proof in [5] of Theorem 1.3 uses an extension
of the ideas behind Theorem 1.1.

\vskip 18pt
\centerline{\caps 2. Proof of the main theorem}
\vskip 6pt

The following appears as Lemmas 2.1 and 2.2
in [4].

\noindent {\bf Lemma 2.1.} {\sl Let $X$ be a Cayley graph
of a finite group, and let $D$ be a non-empty subset
of $V(X)$ such that $|\partial D|/|D| = h(X)$ and
$|D| \leq |V(X)|/2$. Then
$|D| > |V(X)|/4$. Furthermore, the subgraphs induced
by $D$ and its complement $D^c$ are connected.}

\noindent {\sl Proof of Theorem 1.1.} Suppose that
(2) and (3) of 1.1 do not hold. Our aim is to show that
(1) must be true. 

We fix $\epsilon$ to be some real number strictly
between $0$ and ${2 \over \sqrt 3} -1$, but where we
view it as very small.
Since the rank gradient of $(G, \{ G_i \})$ is non-zero,
there is a subgroup $H$ in the lattice $\{ G_i \}$ such that
$(d(H)-1)/[G:H]$ 
is at most $(1 + \epsilon)$ times the rank gradient of $(G, \{ G_i \})$.
The pair $(H, \{ G_i \cap H \})$ has rank gradient
at least $[G:H]$ times the rank
gradient of $(G, \{G_i\})$, since
$\{ G_i \cap H \}$ is a sublattice of $\{ G_i \}$. 
So, $d(H) - 1$ is at most $(1 + \epsilon)$ times the
rank gradient of $(H, \{ G_i \cap H \})$.
Also, (2) does not
hold for this sublattice. Hence, by replacing $G$ by $H$,
and replacing $\{ G_i \}$ by $\{ G_i \cap H \}$,
we may assume that $d(G)-1$ is at most $(1 + \epsilon)$ times the
rank gradient of $(G, \{ G_i \})$.

Let $K$ be a finite 2-complex having 
fundamental group $G$, arising from a minimal
generator presentation of $G$. Thus, $K$ has
a single vertex and $d(G)$ edges. Let $L$ be
the sum of the lengths of the relations in this
presentation. 
Let $K_i \rightarrow K$ be the covering 
corresponding to $G_i$, and let $X_i$ be the
1-skeleton of $K_i$. Since we are assuming that
$G$ does not have Property $(\tau)$ with respect
to $\{ G_i \}$, we may pass to a sublattice where $h(X_i) \rightarrow 0$.
For each $i$, let $D_i$ be a non-empty subset
of $V(X_i)$ such that $|\partial D_i|/|D_i| = h(X_i)$ and
$|D_i| \leq |V(X_i)|/2$. Lemma 2.1
asserts that $|D_i| > |V(X_i)|/4$. We will use $D_i$
to construct a decomposition of $K_i$ into two
overlapping subsets. Let $A_i$ (respectively, $B_i$) 
be the closure of the union of those cells in $K_i$ 
that intersect $D_i$ (respectively, $D_i^c$).
Let $C_i$ be $A_i \cap B_i$.
Lemma 2.1 asserts that the subgraphs induced by $D_i$
and $D_i^c$ are connected, and hence so are $A_i$ and
$B_i$. The 1-skeleton of $A_i$ consists of three
types of edges (that are not mutually exclusive):
\item{(i)} those edges with both endpoints in $D_i$,
\item{(ii)} the edges in $\partial D_i$,
\item{(iii)} those edges in the boundary of a 2-cell that
intersects both $D_i$ and $D_i^c$.

\noindent If we consider the $d(G)$ oriented edges of $K_i$
emanating from the identity vertex in $K_i$ and
translate these edges by the covering translations
in $D_i$ (where we view $D_i$ as subset of $G/G_i$), 
we will cover every edge in (i),
and possibly others. Hence, there are at most
$|D_i| d(G)$ edges of type (i). Similarly,
any type (iii) edge lies in a 2-cell that
intersects both $D_i$ and $D_i^c$. Place the
basepoint of this 2-cell at one its vertices
in $D_i$ that is adjacent to $\partial D_i$.
Translating this basepoint to the identity vertex,
we obtain a corner of a 2-cell incident
to the identity vertex. There are at most
$L$ such corners. The 2-cell runs over at most $L$ 1-cells,
and there are at most $|\partial D_i|$ possibilities
for the translation.
So, there are no more than  $|\partial D_i|L^2$
type (iii) edges. There are $|\partial D_i|$
type (ii) edges, and so, there are
at most $|\partial D_i| (L^2 + 1)$ type (ii) and (iii) edges
in total.
Since $d(\pi_1A_i) -1$ is at most the number
of edges of $A_i$ minus the number of its vertices,
$$\eqalign{
d(\pi_1A_i) -1 
&\leq |D_i| d(G) + |\partial D_i| (L^2+1) -|D_i| \cr 
&= |D_i| (d(G) -1 + h(X_i) (L^2+1)) \cr
&\leq \textstyle{1 \over 2} [G:G_i] (d(G) -1 + h(X_i) (L^2+1)) \cr
&\leq \textstyle{1 \over 2} (1 + \epsilon) [G:G_i] (d(G) -1)
\hbox{ when }h(X_i) \hbox{ is sufficiently small} \cr
&\leq \textstyle{1 \over 2} (1 + \epsilon)^2 (d(G_i) - 1).}
$$
A similar inequality holds for $d(\pi_1B_i) - 1$, but
where ${1 \over 2}$ is replaced throughout by ${3 \over 4}$.
We also note that the 1-skeleton of $C_i$ consists
of type (ii) and type (iii) edges (that may also be of
type (i)), and so, for each component
$C_{i,j}$, of $C_i$,
$$d(\pi_1C_{i,j}) \leq |\partial D_i|(L^2+1) = |D_i| h(X_i) (L^2+1)
\leq  \textstyle {1 \over 2} [G:G_i] h(X_i) (L^2+1).$$

If $C_i$ is disconnected, then $G_i$ is an
HNN extension, giving (1). Thus, we
may assume that $C_i$ is connected.
Then, by the Seifert-Van Kampen theorem, $\pi_1 K_i$
($=G_i$) is the pushout of the diagram
$$
\matrix{
\pi_1C_i & \longrightarrow & \pi_1A_i \cr
\Big\downarrow \cr
\pi_1B_i \cr}$$
where the maps are induced by inclusion. These
homomorphisms need not be injective. However, if we
write ${\rm Im}(\pi_1C_i)$ for the image of
$\pi_1C_i$ in $\pi_1K_i$, and so on, then
$\pi_1K_i$ is the pushout of
$$
\matrix{
{\rm Im}(\pi_1C_i) & \longrightarrow & {\rm Im}(\pi_1A_i) \cr
\Big\downarrow \cr
{\rm Im}(\pi_1B_i). \cr}$$
This follows from a straightforward application
of the universal property of pushouts.
The maps in the above diagram are now injections. When $h(X_i)$ is sufficiently
small, neither ${\rm Im}(\pi_1A_i)$ nor
${\rm Im}(\pi_1B_i)$ can be all of $G_i$. This is because
they then have rank at most 
${3 \over 4} (1 + \epsilon)^2 (d(G_i)-1)+1$,
which is less than $d(G_i)$, when $i$
is sufficiently large, by our assumption
that $\epsilon < {2 \over \sqrt 3} - 1$.
Thus, we deduce that $G_i$ is the non-trivial
amalgamated free product of ${\rm Im}(\pi_1A_i)$ and
${\rm Im}(\pi_1B_i)$ along ${\rm Im}(\pi_1C_i)$. $\square$

The argument above gives rather more, in fact, than is
stated in Theorem 1.1. It immediately implies the
following.

\noindent {\bf Addendum 2.2.} {\sl Theorem 1.1 remains
true if (1) is replaced by:
\item{$1'$.} $G_i$ is an amalgamated free product
$P_i \ast_{R_i} Q_i$ or HNN extension $P_i \ast_{R_i}$
for infinitely many $i$, and in some subsequence,
$$\lim_{i \rightarrow \infty} {d(R_i) \over d(G_i)} = 0.$$
Furthermore, in the case when these $G_i$ are amalgamated
free products $P_i \ast_{R_i} Q_i$, 
we may ensure that the following also hold in
this subsequence:
$$\eqalign{
{1 \over 4} &\leq \liminf_i {d(P_i) \over d(G_i)} \leq 
\limsup_i {d(P_i) \over d(G_i)} \leq
{3 \over 4} \cr
{1 \over 4} &\leq \liminf_i {d(Q_i) \over d(G_i)}  \leq
\limsup_i {d(Q_i) \over d(G_i)} \leq {3 \over 4}. \cr}$$}

\vskip 18pt
\centerline {\caps 3. Rank gradient}
\vskip 6pt

The first examples one comes to of groups with non-zero
rank gradient are free non-abelian groups. 
If $G_i$ is a finite index subgroup of a finitely generated
free group $F$, then
$$d(G_i) - 1 = (d(F) - 1)[F:G_i],$$
and so the rank gradient of $F$ is $d(F)-1$.

Since ${\rm SL}(2, {\Bbb Z})$ has a free non-abelian normal
subgroup of finite index, the following lemma
implies that it has non-zero rank gradient.

\noindent {\bf Lemma 3.1.} {\sl Let $H$ be a finite index normal
subgroup of a finitely generated infinite group $G$, and let $\{ G_i \}$ be a
collection of finite index normal subgroups of
$G$. Then $(G, \{ G_i \})$ has non-zero rank gradient 
if and only if $(H, \{ G_i \cap H \})$ has non-zero rank
gradient. Hence, $G$ has non-zero rank gradient
if and only if the same is true of $H$.}

\noindent {\sl Proof.} This is a consequence of the
following inequalities:
$$\eqalign{
{[G:G_i] \over [G:H]} &\leq [H:G_i \cap H] \leq
[G:G_i] \cr
d(G_i) - d(G_i/G_i \cap H) &\leq d(G_i \cap H) \leq
(d(G_i) - 1)[G_i:G_i \cap H] + 1.}$$
We note that $G_i/G_i \cap H$ is a finite 
group with order at most $[G:H]$,
and hence $d(G_i/G_i \cap H)$ is bounded,
independent of $i$. $\square$

The same lemma gives the following more general
conclusion. Let $G$ be the amalgamated free product
$A \ast_C B$, where $A$ and $B$ are finite, where
$C$ is a proper subgroup of both $A$ and $B$,
and where at least one of $[A:C]$ and $[B:C]$
is more than two. Then $G$ has non-zero rank gradient.

There are two possible generalisations from free non-abelian 
groups. The first
is to groups with deficiency more than one, namely those
groups $G$ admitting a finite presentation $\langle X | R \rangle$
where $|X| > |R|+1$. If we apply the Reidermeister-Schreier
process to a finite index subgroup $G_i$, we obtain a presentation
for $G_i$ with $(|X|-1) [G:G_i] + 1$ generators
and $|R| [G:G_i]$ relations. Hence, the first Betti
number of $G_i$ is at least 
$$(|X| - 1 - |R|)[G:G_i] + 1.$$
This is a lower bound for its rank, and so the
rank gradient of $G$ is at least $(|X| - 1 - |R|)$, which
is positive.

The second way to generalise from the example of free
non-abelian groups is to free products of groups.
Here, we have the following result.

\noindent {\bf Proposition 3.2.} {\sl Let $G$ be the
free product of two non-trivial, finitely generated
groups, not both isomorphic to ${\Bbb Z}/2{\Bbb Z}$. Then
$G$ has non-zero rank gradient.}

\noindent {\sl Proof.} Let $G = A \ast B$, where
$|A| > 2$ and $|B| > 1$. This gives a graph of
groups decomposition of $G$. This lifts to
a graph of groups decomposition for any finite
index normal subgroup $G_i$. The vertex groups covering
$A$ are all isomorphic to $A \cap G_i$. There are $[G:G_i] / [A:A \cap G_i]$ such
vertices. A similar statement holds for the vertices
covering $B$. Since the edge group of $G$ is trivial,
so too are all the edge groups of $G_i$, and there
are $[G:G_i]$ of them. Hence, the first Betti number
of the graph for $G_i$ is
$$[G:G_i] \left( 1 - {1 \over [A:A \cap G_i]} - {1 \over [B:B \cap G_i]} \right) + 1.$$
This is a lower bound for the rank of $G_i$.
Thus, the rank gradient of $G$ is positive,
unless one of $[A : A \cap G_i]$ and 
$[B:B \cap G_i]$ is one for infinitely many $i$ {\sl or} 
they are both two for infinitely many $i$.
In the former case, the vertex groups of $G_i$ that
cover the $A$ (or $B$) group are themselves isomorphic to
$A$ (or $B$), and there are $[G:G_i]$ of them.
In this case, $G_i$ is a free product, with at
least $[G:G_i]$ summands isomorphic to $A$ (or $B$).
By Grushko's theorem [9], the rank of $G_i$ is then
at least $[G:G_i] d(A)$ (or $[G:G_i] d(B)$),
and so $G$ has non-zero rank gradient.
Now consider the second case, where $[A : A \cap G_i]$ and 
$[B:B \cap G_i]$ are both two. Since $A$ has
more than two elements, $A \cap G_i$ is non-trivial.
The vertex groups of $G_i$ covering $A$ are all
isomorphic to $A \cap G_i$, and there are $[G:G_i]/2$
of them. Hence, $d(G_i) \geq d(A \cap G_i) [G:G_i]/2$,
and again $G$ has non-zero rank gradient. $\square$

The condition in the above theorem that $G$ not
be isomorphic to ${\Bbb Z}/2{\Bbb Z} \ast {\Bbb Z}/2{\Bbb Z}$ 
is clearly necessary.
This is because ${\Bbb Z}/2{\Bbb Z} \ast {\Bbb Z}/2{\Bbb Z}$ 
contains ${\Bbb Z}$ as a finite
index normal subgroup, and hence has zero rank gradient,
by Lemma 3.1

There is another source of examples of groups
having collections of subgroups with non-zero
rank gradient.

\noindent {\bf Proposition 3.3.} {\sl Let $G$ be
a finitely generated group that admits a surjective homomorphism
$\phi \colon G \rightarrow F$ onto a free
non-abelian group. Then, for any collection 
$\{ G_i \}$ of finite index subgroups that
each contain the kernel of $\phi$, $(G, \{ G_i \})$ 
has non-zero rank gradient.}

\noindent {\sl Proof.} For any such $G_i$,
$[G:G_i] = [F: \phi(G_i)]$. Hence,
$$d(G_i) \geq d(\phi(G_i)) = (d(F) - 1) [F:\phi(G_i)] +1
= (d(F) - 1)[G:G_i] + 1.$$
So, the rank gradient of $(G, \{ G_i \})$ is at least
$d(F) - 1$. $\square$

This is relevant in 3-manifold theory. 

\noindent {\bf Corollary 3.4.} {\sl Let $M$ be a
compact irreducible 3-manifold with non-empty
boundary, other than an $I$-bundle over a disc, annulus, torus
or Klein bottle. Then $\pi_1 M$
has an infinite lattice of finite index
subgroups with non-zero rank gradient.
Hence, the corresponding covers have non-zero
Heegaard gradient.}

\noindent {\sl Proof.} $\pi_1M$ has a finite
index normal subgroup that admits a surjective homomorphism
$\phi$ onto ${\Bbb Z} \ast {\Bbb Z}$, 
by a theorem of Cooper, Long and Reid [1]. 
The finite index subgroups of $\pi_1M$ that contain
the kernel of $\phi$ form the required infinite
lattice. The final statement of the corollary follows
from the observation that the Heegaard
genus of a 3-manifold is at least the rank
of its fundamental group. $\square$

We now turn to groups with zero rank gradient.
The first collection of examples are mapping tori.
These are constructed from a finitely generated group $A$ and a
homomorphism $\phi \colon A \rightarrow A$. The
associated {\sl mapping torus} $G$ is
$$\langle A, t \ | \ a = t^{-1}\phi(a)t \hbox{ for all } a \in A \rangle.$$
Its rank is at most $d(A) + 1$. It
admits a surjective homomorphism $G \rightarrow {\Bbb Z}$,
sending $A$ to $0$, and $t$ to $1$. Compose
this with the homomorphism to ${\Bbb Z}/n{\Bbb Z}$,
reducing the integers modulo $n$. The kernel of this
homomorphism is a group $G_n$, which is an index
$n$ normal subgroup of $G$. It is isomorphic to the
mapping torus associated with $A$ and $\phi^n$,
and hence has rank at most $d(A) + 1$. Thus,
the collection $\{ G_n \}$ has bounded rank,
and hence $(G, \{ G_n \})$ has zero rank gradient.

Our second class of groups with zero rank gradient
is the $S$-arithmetic groups with trivial
congruence kernel. It is a result
of Sury and Venkataramana [13] that, for 
such a group $G$, there is uniform bound on the rank of its
principal congruence subgroups. In [13],
they proved this in the case of ${\rm SL}(n, {\Bbb Z})$, 
where $n \geq 3$, but they state that the
proof carries over to this much larger
class of groups. Hence, $G$ has zero rank gradient.
I am grateful to the referee for suggesting that we consider
lattices other than ${\rm SL}(n, {\Bbb Z})$ (with
$n > 2$).
An interesting further collection of examples
comes from the following result.

\noindent {\bf Theorem 1.2.} {\sl Any finitely presented, residually finite,
amenable group has non-positive rank gradient.}

\noindent {\sl Proof.} Let $G$ be such a group, which we
may assume is infinite.
Let $\{ G_i \}$ be its finite index normal subgroups.
Since $G$ is infinite, amenable and residually finite, it 
does not have Property $(\tau)$, and so (2)
of Theorem 1.1 cannot hold. Suppose that ($1'$)
of Addendum 2.2 holds. Then, for infinitely many $i$, $G_i$ is a graph of
groups and hence acts cocompactly on a tree.
Each vertex of this tree has stabiliser which
is a conjugate of $P_i$ or $Q_i$. The number
of edges emanating from this vertex is the
index of $R_i$ in $P_i$ or $Q_i$, as appropriate.
By ($1'$) this index can be arbitrarily large,
and so the tree is not homeomorphic to the
real line. Hence, $G$ contains a non-abelian free group,
which contradicts the assumption that it
is amenable.
We deduce that (3) holds: $G$ has zero
rank gradient. $\square$

We close this section with examples of groups,
each having an infinite lattice of subgroups with
non-zero rank gradient, and another infinite
lattice of subgroups with zero rank gradient.
Such groups arise as the fundamental group of
a hyperbolic 3-manifold that fibres over the circle,
with a fibre a compact surface with negative
Euler characteristic and non-empty boundary.
Since these groups are mapping tori, they
have zero rank gradient. But Corollary 3.4 also
gives that they have non-zero rank gradient
with respect to some infinite lattice of finite index 
subgroups. A concrete example is the fundamental group
of the figure-eight knot complement [14]. This is
an index 24 subgroup of ${\rm SL}(2, {\cal O}_3)$,
where ${\cal O}_3$ is the ring of integers
in ${\Bbb Q}(\sqrt{-3})$.

\vskip 18pt
\centerline {\caps 4. Examples}
\vskip 6pt

In this section, we investigate which possible
combinations of (1), (2) and (3) of Theorem 1.1
can arise.

\noindent {\bf Example 4.1.} ${\Bbb Z}$ is trivially
an HNN extension. It does not have Property $(\tau)$.
And it has zero rank gradient. Thus it satisfies
(1) and (3) but not (2). More generally, any
mapping torus has these properties. 

\noindent {\bf Example 4.2.} A non-abelian free
group has non-zero rank gradient. It admits
a surjective homomorphism onto ${\Bbb Z}$
and hence does not have Property $(\tau)$.
Any finite index subgroup is free, and
is, in particular, a non-trivial graph of groups.
So, it satisfies (1), but not (2) or (3).

\noindent {\bf Example 4.3.} ${\rm SL}(2, {\Bbb Z})$ is
an amalgamated free product 
${\Bbb Z}/4{\Bbb Z} \ast_{{\Bbb Z}/2{\Bbb Z}} {\Bbb Z}/6{\Bbb Z}$,
and hence any finite index subgroup is a non-trivial
graph of groups.
It has Property $(\tau)$ with respect to its
congruence subgroups [7]. And we have already seen
that it has non-zero rank gradient. Thus, the
lattice of congruence subgroups
satisfies (1) and (2) but not (3).

\noindent {\bf Example 4.4.} ${\rm SL}(n, {\Bbb Z})$, where
$n \geq 3$, has Property (T). Hence, no finite index
subgroup is either an HNN extension
or an amalgamated free product. Another consequence
is that it has Property $(\tau)$. We have already
seen that it has zero rank gradient. So, these
groups satisfy (2) and (3) but not (1).

\noindent {\bf Example 4.5.} ${\rm SL}(2, {\Bbb Z}[{1 \over p}])$,
where $p$ is a prime,
satisfies (1), (2) and (3). It was proved
by Serre [12] that ${\rm SL}(2, {\Bbb Z}[{1 \over p}])$
decomposes as a non-trivial graph of groups,
and hence, so does any finite index subgroup.
It has infinitely many finite index subgroups,
for example, its principal congruence subgroups,
and hence it satisfies (1). In fact, any
finite index subgroup is contained in a
principal congruence subgroup [11]. In other
words, it has trivial congruence kernel, 
and so, by the remarks about $S$-arithmetic
groups in \S3, it satisfies (3). Also
using this control on its finite index
subgroups, an argument of Lubotzky in [7]
gives that it has Property $(\tau)$,
establishing (2). I am grateful to the
referee for suggesting this example.

\vfill\eject
\noindent {\bf Example 4.6.} The Grigorchuk group $\Gamma$
is finitely generated, infinite, torsion, amenable, and
residually finite [3]. It therefore satisfies
neither (1) nor (2). We claim also that $\Gamma$
has zero rank gradient, and so that (3) is true. 
For each positive integer $k$, $\Gamma$
has a `level $k$ congruence subgroup'
$St_\Gamma(k)$, which is normal and has finite
index [2]. For any $k \geq 4$, $St_\Gamma(k)$
is isomorphic to the product of $2^{k-3}$
copies of $St_\Gamma(3)$. (This follows from
Proposition 30 (iii) and (vi) in [2]). So, $d(St_\Gamma(k))
\leq 2^{k-3} d(St_\Gamma(3))$. However,
the index of $St_\Gamma(k)$ in $\Gamma$
is $2^{5 \cdot 2^{k-3}+2}$. Thus, these
subgroups have zero rank gradient.

However, $\Gamma$ is not a finitely presented example. There are
no known examples of a finitely presented group
that fails to have Property $(\tau)$ but
that does not have a finite index subgroup
with infinite abelianisation. As a result,
there are no known finitely presented groups 
where neither (1) nor (2) hold. Whether such 
groups can exist is an important unresolved
problem.

There is one other theoretical possibility:
where (2) holds but (1) and (3) do not. The 
difficulty here is the problem of establishing that a
group has non-zero rank gradient if it
is not an amalgamated free product. It seems
likely that the absence of known examples here
is merely due to a lack of mathematical
tools, rather than due to genuine non-existence.

\vskip 18pt
\centerline {\caps 5. Finite index subgroups having free
non-abelian quotients}
\vskip 6pt

In this section, we prove the following result,
which is a weaker form of Theorem 1.3.

\noindent {\bf Theorem 5.1.} {\sl Let $G$ be a finitely
presented group. Suppose that there exists a sequence 
$G_1 \geq G_2 \geq \dots$
of finite index subgroups of $G$, each normal in $G_1$, such that
\itemitem{(i)} $G_i/G_{i+1}$ is abelian for all $i \geq 1$;
\itemitem{(ii)} $\lim_{i \rightarrow \infty} 
((\log [G_i : G_{i+1}]) / [G:G_i]) = \infty$;
\itemitem{(iii)} $\limsup_i (d(G_i/G_{i+1}) / [G:G_i])  > 0$.

\noindent Then $G_i$ is a non-trivial graph of
groups for all sufficiently large $i$.
}

We first replace $G$ by $G_1$, so that each $G_i$
is normal in $G$.
Since $G_i / G_{i+1}$ is abelian for all $i \geq 1$, 
the following theorem of Lubotzky and Weiss [8] applies.
This is, in fact, not exactly how they stated their
result (which appears as Theorem 3.6 of [8]), but this
formulation can readily be deduced from their
argument.

\noindent {\bf Theorem.} {\sl Suppose that a finitely
generated group $G$ has Property $(\tau)$ with respect
to a collection $\{ G_i \}$ of finite index subgroups.
Then there is a constant $c$ with the following
property. If $G_{i+1} \triangleleft G_i$ and
$G_i / G_{i+1}$ is abelian, then $|G_i / G_{i+1}| < c^{[G:G_i]}$.}

Hence, by properties (i) and (ii), $G$ does not
have Property $(\tau)$ with respect to $\{ G_i \}$.
Since $\{ G_i \}$ is a nested sequence,
$(d(G_i) - 1)/[G:G_i]$ is a non-increasing
function of $i$. So, 
the rank gradient of $(G, \{ G_i \})$ is
$$\inf_i \left \{ {d(G_i) - 1 \over [G:G_i]} \right \}
= \limsup_i \left \{ {d(G_i) - 1 \over [G:G_i]} \right \}
\geq \limsup_i \left \{ {d(G_i/G_{i+1}) \over [G:G_i]} \right \},$$
which by (iii) is positive.
So, the only possible conclusion of Theorem 1.1
is (1). Once we know that one $G_i$ is a non-trivial
graph of groups, the same is true of all its finite
index subgroups. This proves Theorem 5.1. $\square$

\vskip 18pt
\centerline{\caps References}
\vskip 6pt

\item{1.} {\caps D. Cooper, D. Long, A. Reid,} 
{\sl Essential closed surfaces in bounded $3$-manifolds.} 
J. Amer. Math. Soc. 10 (1997) 553--563.

\item{2.} {\caps P. de la Harpe}, {\sl
Topics in Geometric Group Theory},
Chicago Lect. Math (2000).

\item{3.} {\caps R. Grigorchuk}, {\sl Burnside's
problem on periodic groups}, Functional Anal. Appl.
{14 (1980) 41--43.

\item{4.} {\caps M. Lackenby}, {\sl Heegaard splittings,
the virtually Haken conjecture and Property ($\tau$)},
Preprint.

\item{5.} {\caps M. Lackenby}, {\sl A characterisation
of large finitely presented groups}, Preprint.

\item{6.} {\caps A. Lubotzky}, {\sl Dimension
function for discrete groups}, Proceedings of Groups,
St. Andrews 1985, 254--262, 
London Math. Soc. Lecture Note Ser., 121, 
(CUP, 1986). 

\item{7.} {\caps A. Lubotzky}, {\sl Discrete Groups,
Expanding Graphs and Invariant Measures}, Progr.
in Math. 125 (1994)

\item{8.} {\caps A. Lubotzky, B. Weiss}, {\sl Groups and expanders},
Expanding graphs (Princeton, 1992) 95--109, DIMACS Ser.
Discrete Math. Theoret. Comput. Sci, 10, Amer. Math.
Soc., Providence, RI, 1993.

\item{9.} {\caps R. Lyndon, P. Schupp,} {\sl 
Combinatorial group theory}, Springer-Verlag, Berlin, 1977. 

\item{10.} {\caps G. Margulis}, {\sl Explicit construction
of expanders}, Problemy Peredav di Informacii 9 (1973) 71--80.

\item{11.} {\caps J-P. Serre}, {\sl Le probl\`eme 
des groupes de congruence pour ${\rm SL}_2$}, Ann. Math. 92 (1970) 489--527.

\item{12.} {\caps J-P. Serre}, {\sl Arbres, amalgames, ${\rm SL}_{2}$.}
Ast\'erisque 46 (1977).

\item{13.} {\caps B. Sury and T.N. Venkataramana},
{\sl Generators for all principal congruence subgroups
of ${\rm SL}(n, {\Bbb Z})$ with $n \geq 3$}, Proc. Amer. Math.
Soc. 122 (1994) 355--358.

\item{14.} {\caps W. Thurston}, {\sl The geometry and topology
of 3-manifolds}, Lecture Notes, Princeton, 1978.

\vskip 12pt
\+ Mathematical Institute, Oxford University, \cr
\+ 24-29 St Giles', Oxford OX1 3LB, UK. \cr

\end